\documentclass[graybox]{svmult}

\usepackage{mathptmx}       
\usepackage{helvet}         
\usepackage{courier}        
\usepackage{type1cm}        
\usepackage{makeidx}         
\usepackage{graphicx}        
\usepackage{multicol}        
\usepackage[bottom]{footmisc}

\usepackage{amssymb}
\usepackage{amsfonts}
\usepackage{comment}
\usepackage{mathtools}

\makeindex             

\newcommand{\B}{\mathcal{B}}
\newcommand{\C}{\mathbb{C}}
\newcommand{\K}{\mathbb{K}}
\newcommand{\N}{\mathbb{N}}
\newcommand{\g}{\mathfrak{g}}
\newcommand{\h}{\mathfrak{h}}

\newcommand{\dfn}[1]{\begin{definition} #1 \end{definition}}
\newcommand{\thm}[1]{\begin{theorem} #1 \end{theorem}}
\newcommand{\prop}[1]{\begin{proposition} #1 \end{proposition}}

\newcommand{\prf}[1]{\begin{proof} \smartqed #1 \qed \end{proof}}

\begin{document}

\title*{Classification of Low Dimensional 3-Lie Superalgebras}
\author{Viktor Abramov, Priit L\"att}
\institute{%
  Viktor Abramov \at
  Institute of Mathematics, University of Tartu, \\
  Liivi 2--602, Tartu 50409, Estonia, \\
  \email{viktor.abramov@ut.ee}
  \and
  Priit L\"att \at
  Institute of Mathematics, University of Tartu, \\
  Liivi 2--602, Tartu 50409, Estonia, \\
  \email{priit\_lt@ut.ee}
}
%
%
\maketitle

\abstract*{A notion of $n$-Lie algebra introduced by V.T. Filippov can be viewed as a generalization of a concept of binary Lie algebra to the algebras with $n$-ary multiplication law. A notion of Lie algebra can be extended to $\mathbb Z_2$-graded structures giving a notion of Lie superalgebra. Analogously a notion of $n$-Lie algebra can be extended to $\mathbb Z_2$-graded structures by means of a graded Filippov identity giving a notion of $n$-Lie superalgebra. We propose a classification of low dimensional 3-Lie superalgebras. We show that given an $n$-Lie superalgebra equipped with a supertrace one can construct the $(n+1)$-Lie superalgebra which is referred to as the induced $(n+1)$-Lie superalgebra. A Clifford algebra endowed with a $\mathbb Z_2$-graded structure and a graded commutator can be viewed as the Lie superalgebra. It is well known that this Lie superalgebra has a matrix representation which allows to introduce a supertrace. We apply the method of induced Lie superalgebras to a Clifford algebra to construct the 3-Lie superalgebras and give their explicit description by ternary commutators.}

\abstract{A notion of $n$-Lie algebra introduced by V.T. Filippov can be viewed as a generalization of a concept of binary Lie algebra to the algebras with $n$-ary multiplication law. A notion of Lie algebra can be extended to $\mathbb Z_2$-graded structures giving a notion of Lie superalgebra. Analogously a notion of $n$-Lie algebra can be extended to $\mathbb Z_2$-graded structures by means of a graded Filippov identity giving a notion of $n$-Lie superalgebra. We propose a classification of low dimensional 3-Lie superalgebras. We show that given an $n$-Lie superalgebra equipped with a supertrace one can construct the $(n+1)$-Lie superalgebra which is referred to as the induced $(n+1)$-Lie superalgebra. A Clifford algebra endowed with a $\mathbb Z_2$-graded structure and a graded commutator can be viewed as the Lie superalgebra. It is well known that this Lie superalgebra has a matrix representation which allows to introduce a supertrace. We apply the method of induced Lie superalgebras to a Clifford algebra to construct the 3-Lie superalgebras and give their explicit description by ternary commutators.}

\section{Introduction}

Recently, there was markedly increased interest of theoretical physics towards the algebras with $n$-ary multiplication law. Due to the fact that the Lie algebras play a crucial role in theoretical physics, it seems that development of $n$-ary analog of a concept of Lie algebra is especially important. In \cite{Filippov} V.T. Filippov proposed a notion of $n$-Lie algebra which can be considered as a possible generalization of a concept of Lie algebra to structures with $n$-ary multiplication law. In approach proposed by V.T. Filippov an $n$-ary commutator of $n$-Lie algebra is skew-symmetric and satisfies an $n$-ary analog of Jacobi identity which is now called Filippov identity. It is worth to mention that there is an approach different from the one proposed by V.T. Filippov, where a ternary commutator is not skew-symmetric but it obeys a symmetry based on a representation of the group of cyclic permutations $\mathbb Z_3$ by cubic roots of unity \cite{Abramov-Kerner-LeRoy}. It is well known that a concept of Lie algebra can be extended to $\mathbb Z_2$-graded structures with the help of graded commutator and graded Jacoby identity, and a corresponding structure is known under the name of Lie superalgebra.

In the present paper we show that a notion of $n$-Lie algebra proposed by V.T. Filippov can be extended to $\mathbb Z_2$-graded structures by means of graded $n$-commutator and a graded analog of Filippov identity. This $\mathbb Z_2$-graded $n$-Lie algebra will be referred to as a $n$-Lie superalgebra. We show that a method of induced $n$-Lie algebras proposed in \cite{silvestrov2014} and based on an analog of a trace can be applied to $n$-Lie superalgebras if instead of a trace we will be using a supertrace. We introduce the notions such as an ideal of $n$-Lie superalgebra, subalgebra of $n$-Lie superalgebra, descending series and prove several results analogous to the results proved in \cite{silvestrov2014} for $n$-Lie algebras. We propose a classification of low dimensional 3-Lie superalgebras and find their commutation relations. A Clifford algebra can be used to construct a Lie superalgebra if one equips it with a grade commutator. This Lie superalgebra has a matrix representation called supermodule of spinors and this representation can be endowed with a supertrace. Thus we have all basic components of a method of induced $n$-Lie superalgebras and applying this method we construct a series of 3-Lie superalgebras.

\section{Supertrace and Induced $n$-Lie Superalgebras}

A notion of Lie algebra can be extended from binary algebras to algebras with $n$-ary multiplication law with the help of a notion of $n$-Lie algebra, where $n$ is any integer greater or equal to 2. This approach was proposed by V. T. Filippov in \cite{Filippov}, and it is based on $n$-ary analog of Jacobi identity which is now called the Filippov identity. It is well known that a concept of binary Lie algebra can be extended to $\mathbb Z_2$-graded structures giving a notion of Lie superalgebra. Similarly a notion of $n$-Lie algebra can be extended to $\mathbb Z_2$-graded structures giving a structure which we call an $n$-Lie superalgebra \cite{abramov2014,Daletskii-Kushnirevitch}. In this section we give the definitions of $n$-Lie algebra, $n$-Lie superalgebra and show that a structure of induced $n$-Lie algebra based on an analog of trace \cite{silvestrov2014} can be extended to $n$-Lie superalgebras with the help of supertrace.

\dfn{
  Vector space $\g$ endowed with a mapping $[\cdot, \dots, \cdot] : \g^n \to \g$ is said to be a \emph{$n$-Lie algebra}, if $[\cdot, \dots, \cdot]$ is $n$-linear, skew-symmetric and satisfies the identity
  \begin{equation} \label{id:filippov}
    [ x_1, \dots, x_{n-1}, [y_1, \dots, y_n]] = \sum_{i=1}^n [ y_1, \dots, [ x_1, \dots, x_{n-1}, y_i ], \dots, y_n],
  \end{equation}
where $x_1, \dots, x_{n-1}, y_1, \dots, y_n \in \g$.
}
\noindent
In the definition of $n$-Lie algebra the identity (\ref{id:filippov}) is called the Filippov identity \cite{Filippov}. It is clear that for $n = 2$ Filippov identity yields the classical Jacobi identity of a binary Lie algebra.

\dfn{
  Let $\phi : \g^n \to \g$. A linear map $\tau : \g \to \K$ will be referred to as a {$\phi$-trace} if $\tau\left( \phi(x_1, \dots, x_n) \right) = 0$ for all $x_1, \dots, x_n \in \g$.
}
\noindent
In \cite{silvestrov2014} the authors proposed a method based on {$\phi$-trace} which can be applied to an $n$-Lie algebra to construct the $(n+1)$-Lie algebra.

\thm{
  Let $(\g, [\cdot, \dots, \cdot])$ be an $n$-Lie algebra and $\tau$ be a $[\cdot, \dots, \cdot]$-trace. Define $[\cdot, \dots, \cdot]_\tau : \g^{n+1} \to \g$ by
  \begin{equation}
    [x_1, \dots, x_{n+1}]_\tau = \sum_{i=1}^{n+1} (-1)^i \tau(x_i)
    [x_1, \dots, x_{i-1}, x_{i+1}, \dots, x_{n+1}].
  \end{equation}
  Then $(\g, [\cdot, \dots, \cdot]_\tau)$ is the $(n+1)$-Lie algebra.
}
\noindent
It was shown in \cite{abramov2014} that a similar method based on a notion of a supertrace can be used in the case of $n$-Lie superagebras, and given an $n$-Lie superalgebra one can apply this method to induce the $(n+1)$-Lie superalgebra. Let us remind that super vector space is a direct sum of two vector spaces, i.e. $V = V_{\overline{0}} \oplus V_{\overline{1}}$. The dimension of finite dimensional super vector space is denoted as $m|n$ if $dim V_{\overline{0}} = m$ and $dim V_{\overline{0}} = n$. Element $x \in V \setminus \{0\}$ is said to be homogeneous if either $x \in V_{\overline{0}}$ or $x \in V_{\overline{1}}$. For homogeneous elements we can define parity by
\begin{equation}
  |x| = \left\{\begin{array}{r@{}l@{\qquad}l}
    {\overline{0}}, \quad x \in V_{\overline{0}}, \\
    {\overline{1}}, \quad x \in V_{\overline{1}}.
  \end{array}\right.
\end{equation}
\noindent
In what follows we will assume that element $x$ is homogeneous whenewer $|x|$ is used.

\dfn{
  We say that super vector space $\g$ endowed with $n$-linear map $[\cdot, \dots, \cdot] : \g^n \to \g$ is \emph{$n$-Lie superalgebra} if for all $x_1, \dots, x_n, y_1, \dots, y_{n-1} \in \g$
  \begin{enumerate}
    \item $|[x_1, \dots, x_n]| = \sum_{i=1}^{n} |x_i|$, \vspace{0.1cm}
    \item $[x_1, \dots, x_i, x_{i+1}, \dots, x_n] = -(-1)^{|x_i||x_{i+1}|} [x_1, \dots, x_{i+1}, x_i, \dots, x_n]$, \vspace{0.1cm}
    \item $[y_1, \dots, y_{n-1}, [x_1, \dots, x_n]] =$ \\
      $= \sum_{i=1}^{n} (-1)^{|\mathrm{x}|_{i-1} |\mathrm{y}|_{n-1}} [x_1, \dots, x_{i-1}, [y_1, \dots, y_{n-1}, x_i], x_{i+1}, \dots, x_n]$,
  \end{enumerate}
  where $\mathrm{x} = (x_1, \dots, x_n)$, $\mathrm{y} = (y_1, \dots, y_{n-1})$,
  and $|\mathrm{x}_i| = \sum_{j=1}^i |x_j|$.
}

\dfn{
  Let $V = V_{\overline{0}} \oplus V_{\overline{1}}$ be a super vector space and let $\phi : V^n \to V$. We say that linear map $S : V \to \K$ is a \emph{$\phi$-supertrace} if
  \begin{enumerate}
    \item $S\left( \phi(x_1, \dots, x_n) \right) = 0$ for all $x_1, \dots, x_n \in V$, \vspace{0.1cm}
    \item $S(x) = 0$ for all $x \in V_{\overline{1}}$.
  \end{enumerate}
}
\noindent
Given an $n$-Lie superalgebra endowed with a supertrace (which satisfies the conditions of the previous definition with respect to a graded commutator of this algebra) we can construct the $(n+1)$-Lie superalgebra by means of a method described in \cite{abramov2014}.

\thm{\label{thm:induced-superalgebra}
  Let $(\g, [\cdot, \dots, \cdot])$ be a $n$-Lie superalgebra and let $S : \g \to \K$ be a $[\cdot, \dots, \cdot]$-supertrace. Define $[\cdot, \dots, \cdot]_S : \g^{n+1} \to \g$ by
  \begin{equation}
    [x_1, \dots, x_{n+1}]_S = \sum_{i=1}^{n+1} (-1)^{i-1} (-1)^{|x_i||\mathrm{x}|_{i-1}} S(x_i) [x_1, \dots, x_{i-1}, x_{i+1}, \dots, x_{n+1}].
    \label{n super-to-n+1 super}
  \end{equation}
  Then $(\g, [\cdot, \dots, \cdot]_S)$ is a $(n+1)$-Lie superalgebra.
}

\section{Properties of Induced $n$-Lie Superalgebras}

In this section we study a structure of induced $n$-Lie superalgebra, and introducing the notions such as ideal of an $n$-Lie superalgebra, derived series, subalgebra of $n$-Lie superalgebra we prove several results which are analogous to the results proved in \cite{silvestrov2014} in the case of $n$-Lie algebras.

\dfn{
  Let $(\g, [\cdot, \dots, \cdot])$ be a $n$-Lie superalgebra and let $\h$ be
  a subspace of $\g$. We say that $\h$ is an \emph{ideal} of $\g$, if for all
  $h \in \h$ and for all $x_1, \dots, x_{n-1} \in \g$, it holds that $[h, x_1, \dots, x_{n-1}] \in \h$.
}

\dfn{
  Let $(\g, [\cdot, \dots, \cdot])$ be a $n$-Lie superalgebra and let $\h$ be an ideal of $\g$. \emph{Derived series} of $\h$ is defined as
  \[ D^{0} (\h) = \h \quad \mbox{and} \quad D^{p+1} (\h) = [D^p (\h), \dots, D^p (\h)], \quad p \in \N, \]
  and the \emph{descending central series} of $\h$ as
  \[ C^{0} (\h) = \h \quad \mbox{and} \quad C^{p+1} (\h) = [C^p (\h), \h, \dots, \h], \quad p \in \N. \]
}
\noindent
An ideal $\h$ of $n$-Lie superalgebra $\g$ is said to be \emph{solvable} if there exists $p \in \N$ such that $D^p (\h) = \{0\}$, and we call $\h$ \emph{nilpotent} if $C^{p} (\h) = \{0\}$ for some $p \in \N$.

\prop{
  Let $(\g, [\cdot, \dots, \cdot])$ be a $n$-Lie superalgebra and let $\h \subset \g$ be a subalgebra. If $S$ is supertrace of $[\cdot, \dots, \cdot]$, then $\h$ is also subalgebra of $(\g, [\cdot, \dots, \cdot]_S)$.
}

\prf{
  Let $\h$ be subalgebra of $n$-Lie superalgebra $(\g, [\cdot, \dots, \cdot])$, $x_1, \dots, x_{n+1} \in \h$ and assume $S$ is a supertrace of $[\cdot, \dots, \cdot]$. Then $[x_1, \dots, x_{n+1}]_S$ is a linear combinaton of elements of $\h$ as desired.
}

\prop{
  Let $\h$ be an ideal of $(\g, [\cdot, \dots, \cdot])$ and assume $S$ is supertrace of $[\cdot, \dots, \cdot]$. Then $\h$ is ideal of $(\g, [\cdot, \dots, \cdot]_S)$ if and only if $[\g, \g, \dots, \g] \subseteq \h$ or $\h \subseteq \ker S$.
}

\prf{
  Let $h \in \h$ and $x_1, \dots, x_n \in \g$. Then
  \begin{eqnarray}
    [x_1, \dots, x_n, h]_S =
    \sum_{i=1}^{n} (-1)^{i-1} (-1)^{|x_i||\mathrm{x}|_{i-1}}
      S(x_i) [x_1, \dots, x_{i-1}, x_{i+1}, \dots, x_n, h] + \\
    (-1)^{n} (-1)^{|h||\mathrm{x}|_{n}} S(h) [x_1, \dots, x_n].
  \end{eqnarray}
\noindent
Since $\h$ is ideal we have $[x_1, \dots, x_{i-1}, x_{i+1}, \dots, x_n, h] \in \h$ for all $i = 1, \dots, n$. Thus $[x_1, \dots, x_n, h]_S \in \h$ is equivalent to
\[ (-1)^{n} (-1)^{|h||\mathrm{x}|_{n}} S(h) [x_1, \dots, x_n] \in \h, \]
 which clearly holds when $S(h) = 0$ or $[x_1, \dots, x_n] \in \h$.
}

\prop{
  Let $(\g, [\cdot, \dots, \cdot])$ be $n$-Lie superalgebra and let $S$ be supertrace of $[\cdot, \dots, \cdot]$. Then induced $(n+1)$-Lie superalgebra $(\g_S, [\cdot, \dots, \cdot]_S)$ is solvable.
}

\prf{
  Assume $(\g, [\cdot, \dots, \cdot])$ is a $n$-Lie superalgebra and $S$ is
  supertrace of $[\cdot, \dots, \cdot]$, and let $x_1, \dots, x_{n+1} \in D^1 (\g_S)$.

  Then for every $i = 1, \dots, n+1$ we have $x_i^1, \dots, x_i^{n+1} \in \g$
  such that $x_i = [x_i^1, \dots, x_i^{n+1}]_S$, in which case
  \begin{eqnarray*}
    [x_1, \dots, x_{n+1}]_S &=& \\
    \sum_{i=1}^{n+1} &(-1)^{i-1} &(-1)^{|x_i||\mathrm{x}|_{i-1}} S([x_i^1, \dots, x_i^{n+1}]_S) [x_1, \dots, x_{i-1}, x_{i+1}, \dots, x_{n+1}] = 0.
  \end{eqnarray*}
}
\noindent
In the light of the last proposition we can immediately see that if $(\g, [\cdot, \dots, \cdot])$ is an $n$-Lie superalgebra, then for the induced $(n+1)$-Lie superalgebra it holds $D^p (\g_S) = \{0\}$, whenever $p \geq 2$.

\prop{
  Let $(\g, [\cdot, \dots, \cdot])$ be $n$-Lie superalgebra,  $S$ supertrace of $[\cdot, \dots, \cdot]$ and assume $(n+1)$-Lie superalgebra $(\g_S, [\cdot, \dots, \cdot]_S)$ is induced by $S$. Denote descending central series of $\g$ by $\left(C^p(\g)\right)_{p=0}^\infty$ and denote descending central series of $\g_S$ by $\left(C^p(\g_S)\right)_{p=0}^\infty$. Then
  \[ C^p(\g_S) \subseteq C^p(\g) \quad \mbox{for all} \quad p \in \N. \]
  If there exists $g \in \g$ such that $[g, x_1, \dots, x_n]_S = [x_1, \dots, x_n]$ holds for all $x_1, x_2, \dots, x_n \in \g$, then
  \[ C^p(\g_S) = C^p(\g) \quad \mbox{for all} \quad p \in \N. \]
}

\prf{
  Case $p = 0$ is trivial. Note that for $p = 1$ any $x = [x_1, \dots, x_{n+1}]_S \in C^1 (\g_S)$ can be expressed as
  \[ x = \sum_{i=1}^{n+1} (-1)^{i-1} (-1)^{|x_i||\mathrm{X}|_{i-1}} S(x_i) [x_1, \dots, x_{i-1}, x_{i+1}, \dots, x_{n+1}], \]
  meaning $x$ is a linear combination in $C^1 (\g)$.

  Assume now that there exists $g \in \g$ such that for all $y_1, \dots, y_n \in \g$ it holds that $[g, y_1, \dots, y_n]_S = [y_1, \dots, y_n]$. Then $x = [x_1, \dots, x_n] \in C^1 (\g)$ can be written as $[g, x_1, \dots, x_n]_S$ and thus $x \in C^1 (\g_S)$.

  Next assume that the statement holds for some $p \in \N$ and let $x \in C^{p+1} (\g_S)$. Then there are $x_1, \dots, x_n \in \g$ and $g \in C^p (\g_S)$ such that
  \begin{eqnarray*}
    x = [g, x_1, \dots, x_n]_S = (-1)^{n + |g||\mathrm{x}|_n} [x_1, \dots, x_n, g]_S = \\
    = (-1)^{n + |g||\mathrm{x}|_n} \sum_{i=1}^{n} (-1)^{i-1} (-1)^{|x_i||\mathrm{x}|_{i-1}} S(x_i) [x_1, \dots, x_{i-1}, x_{i+1}, \dots, x_n, g],
  \end{eqnarray*}
  since $g$ can be expressed as a bracket of some elements, and hence $S(g) = 0$. On the other hand, as $g \in C^p (\g_S)$, by our inductive assumption $g \in C^p (\g)$, and thus $x \in C^{p+1} (\g)$.

  To complete the proof, assume that there exists $g \in \g$ such that for all $y_1, \dots, y_n \in \g$ equality $[g, y_1, \dots, y_n]_S = [y_1, \dots, y_n]$ holds. If $x \in C^{p+1} (\g)$, then $x = [h, x_1, \dots, x_{n-1}]$, where $x_1, \dots, x_{n-1} \in \g$ and $h \in C^{p} (\g)$. Altogether we have
  \[ x = [h, x_1, \dots, x_{n-1}] = [g, h, x_1, \dots, x_{n-1}]_S = -(-1)^{|g||h|} [h, g, x_1, \dots, x_{n-1}]_S. \]
  At the same time $h \in C^{p} (\g) = C^p (\g_S)$, which gives us $[h, g, x_1, \dots, x_{n-1}]_S \in C^{p+1} (\g_S)$, meaning $x \in C^{p+1} (\g_S)$, as desired.
}

\section{Low Dimensional Ternary Lie Superalgebras}
In this section we propose a classification of low dimensional ternary Lie superalgebras.
First of all we find the number of different (non-isomorphic) $3$-Lie superalgebras over $\C$ of dimension $m|n$ where $m+n < 5$. We also find the explicit commutation relations of these 3-Lie superalgebras. We use a method which is based on the structure constants of an $n$-Lie superalgebra.

\dfn{
  Let $\g = \g_{\overline{0}} \oplus \g_{\overline{1}}$ be $n$-Lie superalgebra, denote
  \[ \B = \{ e_1, \dots, e_m, f_1, \dots, f_n \} \]
  and assume $\{ e_1, \dots, e_m \}$ spans $\g_{\overline{0}}$ and $\{ f_1, \dots, f_n \}$ spans $\g_{\overline{1}}$. Elements $K_{A_1 \dots A_n}^B$ defined by
  \[ [z_{A_1}, \dots, z_{A_n}] = K_{A_1 \dots A_n}^B z_B, \]
  where $z_{A_1}, \dots, z_{A_n}, z_B \in \B$, are said to be \emph{structure constants} of $\g$ with respect to $\B$.
}

Assume we have a $3$-Lie superalgebra $(\g, [\cdot, \cdot, \cdot])$ of dimension $m|n$ over $\C$. Denote
\[ \B = \{ e_1, \dots, e_m, f_1, \dots, f_n \} = \{ z_1, \dots, z_{m+n} \} \]
and assume $e_\alpha$, $1 \leq \alpha \leq m$, spans the even part of $\g$ and
$f_i$, $1 \leq i \leq n$, spans the odd part of $\g$. Additionaly, let $z_A = e_A$, when $1 \leq A \leq m$, and $z_A = f_{A-m}$, when $m < A \leq m+n$. Since $|[z_1, z_2, z_3]| = |z_1| + |z_2| + |z_3|$ we can express the values of commutator on generators using structure constants in the following form:
\[
  \begin{array}{l}
    \left[e_\alpha, e_\beta, e_\gamma \right] = K_{\alpha \beta \gamma}^\lambda e_\lambda, \\[0.1cm]
    \left[e_\alpha, e_\beta, f_i\right] = K_{\alpha \beta i}^j f_j, \\[0.1cm]
    \, \left[e_\alpha, f_i, f_j\right] = K_{\alpha i j}^\beta e_\beta, \\[0.1cm]
    \, \left[f_i, f_j, f_k\right] = K_{i j k}^l f_l,
  \end{array}
\]
where $\alpha \leq \beta \leq \gamma$ and $i \leq j \leq k$. As all other possible orderings and combinations of generators can be transformed into one of these four forms by graded skew-symmetry of $[\cdot, \cdot, \cdot]$, we will not consider them.

As a next step we can eliminate the combinations that are trivial. To find such brackets we can observe different permutations of arguments. If some permutaion yields the initial ordering without preserving the sign, then this bracket must be zero, as in
\[ [e_1, e_1, f_i] = -(-1)^{|e_1||e_1|} [e_1, e_1, f_1] = - [e_1, e_1, f_i]. \]

Finally we can use the graded Filippov identity. Observe $[z_A, z_B, z_C] = K_{ABC}^D z_D \neq 0$, where $1 \leq A \leq B \leq C \leq m + n$, and calculate
\[ [z_E, z_F, [z_A, z_B, z_C]] \]
using two different paths. Firstly use what is known and write
\[ \left[ z_E, z_F, [z_A, z_B, z_C] \right] = K_{ABC}^D [z_E, z_F, z_D]. \]
Then transform bracket $[z_E, z_F, z_D]$ to $(-1)^{\circlearrowright_{DEF}} [z_{D'}, z_{E'}, z_{F'}]$, where $\{D, E, F\} = \{D', E', F'\}$, but $D' \leq E' \leq F'$, and $(-1)^{\circlearrowright_{DEF}}$ gives the sign that comes from graded skew-symmetry. Note that $[z_{D'}, z_{E'}, z_{F'}]$ can be expressed using structure constants and generators as well and thus we have $[z_{D'}, z_{E'}, z_{F'}] = K_{D' E' F'}^H z_H$, which means that on the one hand
\[ \left[ z_E, z_F, [z_A, z_B, z_C] \right] = (-1)^{\circlearrowright_{DEF}} K_{ABC}^D K_{D' E' F'}^H z_H. \]

On the other hand we can use Filippov identity to calculate $\left[ z_E, z_F, [z_A, z_B, z_C] \right]$:
\begin{eqnarray*}
    \left[ z_E, z_F, [z_A, z_B, z_C] \right] = \left[[z_E, z_F, z_A], z_B, z_C\right] +  (-1)^{|z_A|(|z_E|+|z_F|)} \left[z_A, [z_E, z_F, z_B], z_C \right] + \\
    (-1)^{(|z_A|+|z_B|)(|z_E|+|z_F|)}
        \left[z_A, z_B, [z_E, z_F, z_C]. \right]
\end{eqnarray*}
In every summand we can apply the same construction as described above. To do that, let us denote $z_{AEF} = |z_A|(|z_E|+|z_F|)$ and $z_{ABEF} = (|z_A|+|z_B|)(|z_E|+|z_F|)$. Now reorder the arguments in increasing order and replace the result with linear combination of generators and structure constants. By doing so we end up having
\[
  \begin{array}{l}
    \left[ z_E, z_F, [z_A, z_B, z_C] \right] = (-1)^{\circlearrowright_{AEF} + \circlearrowright_{B'C'G'}} K_{A' E' F'}^G K_{B' C' G'}^H \,z_H + \\[0.1cm] \qquad \qquad
    (-1)^{z_{AEF} + \circlearrowright_{BEF} + \circlearrowright_{A'C'G'}}
        K_{B' E' F'}^G K_{A' C' G'}^H
        \,z_H + \\[0.1cm] \qquad \qquad
    (-1)^{z_{ABEF} + \circlearrowright_{CEF} + \circlearrowright_{A'B'G'}}
        K_{C' E' F'}^G K_{A' B' G'}^H
        \,z_H.
  \end{array}
\]
In other words the following system of quadratic equations emerges:
\[
  \begin{array}{l}
    (-1)^{\circlearrowright_{DEF}} K_{ABC}^D K_{D' E' F'}^H z_H =\
    (-1)^{\circlearrowright_{AEF} + \circlearrowright_{B'C'G'}}
        K_{A' E' F'}^G K_{B' C' G'}^H
        \,z_H + \\[0.1cm] \qquad \qquad
    (-1)^{z_{AEF} + \circlearrowright_{BEF} + \circlearrowright_{A'C'G'}}
        K_{B' E' F'}^G K_{A' C' G'}^H
        \,z_H + \\[0.1cm] \qquad \qquad
    (-1)^{z_{ABEF} + \circlearrowright_{CEF} + \circlearrowright_{A'B'G'}}
        K_{C' E' F'}^G K_{A' B' G'}^H
        \,z_H,
  \end{array}
\]
where generators $z_H$ are known and structure constants $K_{ABC}^{D}$ are unknown. Further more, for every $H \in \{1, 2, \dots, m+n \}$ we have
\[
  \begin{array}{l}
    (-1)^{\circlearrowright_{DEF}} K_{ABC}^D K_{D' E' F'}^H =\
    (-1)^{\circlearrowright_{AEF} + \circlearrowright_{B'C'G'}} K_{A' E' F'}^G K_{B' C' G'}^H + \\[0.1cm] \qquad \qquad
    (-1)^{z_{AEF} + \circlearrowright_{BEF} + \circlearrowright_{A'C'G'}}
        K_{B' E' F'}^G K_{A' C' G'}^H + \\[0.1cm] \qquad \qquad
    (-1)^{z_{ABEF} + \circlearrowright_{CEF} + \circlearrowright_{A'B'G'}}
        K_{C' E' F'}^G K_{A' B' G'}^H,
  \end{array}
\]

In summary, we have a system of quadratic equations whose solutions are possible structure constants for $m|n$-dimensional $3$-Lie superalgebra. We note however, that the structure constants are depending on the choice of basis for the super vector space and thus invariant solutions have to be removed case by case.

Applying the described algorithm to concrete cases gives us the following theorems.

\thm{
  $3$-Lie superalgebras over $\C$, whose super vector space dimension is $0|1$ or $1|1$, is Abelian.
}

\thm{
  $3$-Lie superalgebras over $\C$, whose super vector space dimension is
  $0|2$ or $1|2$, are either Abelian or isomorphic to $3$-Lie superalgebra
  $\h$ whose non-trivial commutation relations are
  \[
    \left\{
      \begin{array}{l}
        \left[f_1, f_1, f_1\right] = -f_1 + f_2, \\
        \left[f_1, f_1, f_2\right] = -f_1 + f_2, \\
        \left[f_1, f_2, f_2\right] = -f_1 + f_2, \\
        \left[f_2, f_2, f_2\right] = -f_1 + f_2, \\
      \end{array}
    \right.
    \quad \mbox{or} \quad
    \left[f_1, f_1, f_1\right] = f_2,
  \]
  where $f_1, f_2$ are odd generators of $\h$.
}

\thm{
  $3$-Lie superalgebras over $\C$, whose super vector space dimension is
  $2|1$, are either Abelian or isomorphic to $3$-Lie superalgebra
  $\h$ whose non-trivial commutation relations are
  \[
    \left\{
      \begin{array}{l}
        \left[e_1, f_1, f_1\right] = e_1 + e_2, \\
        \left[e_2, f_1, f_1\right] = -e_1 - e_2, \\
      \end{array}
    \right.
    \quad
    [e_1, e_1, f_1] = f_1,
    \quad \mbox{or} \quad
    [f_1, f_1, f_1] = f_1,
  \]
  where $e_1, e_2$ are even generators of $\h$ and $f_1$ is odd generator
  of $\h$.
}
\section{Supermodule over Clifford algebra}
In this section we apply the method described in Section 2 to a Clifford algebra. It is well known that a Clifford algebra can be equipped with the structure of superalgebra if one associates degree 1 to each generator of Clifford algebra and defines the degree of product of generators as the sum of degrees of its factors. Then making use of a graded commutator we can consider a Clifford algebra as the Lie superalgebra. A Clifford algebra has a matrix representation and this allows to introduce a supertrace. Hence we have a Lie superalgebra endowed with a supertrace, and we can apply the method described in Theorem 2 to construct the 3-Lie superalgebra. In this section we will give an explicit description of the structure of this constructed 3-Lie superalgebra.

A Clifford algebra $C_n$ is the unital associative algebra over $\mathbb C$ generated by $\gamma_1,\gamma_2,\ldots,\gamma_n$ which obey the relations
\begin{equation}
\gamma_i\gamma_j+\gamma_j\gamma_i=2\,\delta_{ij}\,e,\quad i,j=1,2,\ldots,n,
\end{equation}
where $e$ is the unit element of Clifford algebra. Let ${\mathcal N}=\{1,2,\ldots,n\}$ be the set of integers from 1 to $n$. If $I$ is a subset of $\mathcal N$, i.e. $I=\{i_1,i_2,\ldots,i_k\}$ where $1\leq i_1<i_2<\ldots<i_k\leq n$, then one can associate to this subset $I$ the monomial $\gamma_I=\gamma_{i_1}\gamma_{i_2}\ldots\gamma_{i_k}$. If $I=\emptyset$ one defines $\gamma_{\emptyset}=e$. The number of elements of a subset $I$ will be denoted by $|I|$. It is obvious that the vector space of Clifford algebra $C_n$ is spanned by the monomials $\gamma_I$, where $I\subseteq {\mathcal N}$. Hence the dimension of this vector space is $2^n$ and any element $x\in C_n$ can be expressed in terms of these monomials as
$$
x=\sum_{I\subseteq {\mathcal N}} a_I\gamma_I,
$$
where $a_I=a_{i_1i_2\ldots i_k}$ is a complex number. It is easy to see that one can endow a Clifford algebra $C_n$ with the ${\mathbb Z}_2$-graded structure by assigning the degree $|\gamma_I|=|I|\,(\mbox{mod}\,2)$ to monomial $\gamma_I$. Then a Clifford algebra $C_n$ can be considered as the superalgebra since for any two monomials it holds $|\gamma_I\gamma_J|=|\gamma_I|+|\gamma_J|$.

Another way to construct this superalgebra which does not contain explicit reference to Clifford algebra is given by the following theorem.
\begin{theorem}
Let $I$ be a subset of ${\mathcal N}=\{1,2,\ldots,n\}$, and $\gamma_I$ be a symbol associated to $I$. Let $C_n$ be the vector space spanned by the symbols $\gamma_I$. Define the degree of $\gamma_I$ by $|\gamma_I|=|I| (\,\mbox{mod}\,2)$, where $|I|$ is the number of elements of $I$, and the product of $\gamma_I,\gamma_J$ by
\begin{equation}
\gamma_I\,\gamma_J=(-1)^{\sigma(I,J)}\gamma_{I\Delta J},
\label{superalgebra product}
\end{equation}
where $\sigma(I,J)=\sum_{j\in J}\sigma(I,j)$, $\sigma(I,j)$ is the number of elements of $I$ which are greater than $j\in J$, and $I\Delta J$ is the symmetric difference of two subsets. Then $C_n$ is the unital associative superalgebra, where the unit element $e$ is $\gamma_\emptyset$.
\end{theorem}
This theorem can be proved by means of the properties of symmetric difference of two subsets. We remind a reader that the symmetric difference is commutative $I\oplus J=J\oplus I$, associative $(I\Delta J)\Delta K=I\Delta (J\Delta K)$ and $I\Delta \emptyset=\emptyset\Delta I$. The latter shows that $\gamma_{\emptyset}$ is the unit element of this superalgebra. The symmetric difference also satisfies $|I\Delta J|=|I|+|J|\;(\mbox{mod}\,2)$. Hence $C_n$ is the superalgebra.

The superalgebra $C_n$ can be considered as the super Lie algebra if for any two homogeneous elements $x,y$ of this superalgebra one introduces the graded commutator $[x,y]=xy-(-1)^{|x||y|}yx$ and extends it by linearity to a whole superalgebra $C_n$. We will denote this super Lie algebra by ${\frak C}_n$. Then $\{\gamma_I\}_{I\subseteq {\mathcal N}}$ are the generators of this super Lie algebra ${\frak C}_n$, and its structure is entirely determined by the graded commutators of $\gamma_I$. Then for any two generators $\gamma_I,\gamma_J$ we have
\begin{equation}
[\gamma_I,\gamma_J]=f(I,J)\;\gamma_{I\Delta J},
\label{binary commutators}
\end{equation}
where $f(I,J)$ is the integer-valued function of two subsets of $\mathcal N$ defined by
$$
f(I,J)=(-1)^{\sigma(I,J)}\big(1-(-1)^{|I\cap J|}\big),
$$
It is easy to verify that the degree of graded commutator is consistent with the degrees of generators, i.e. $[\gamma_I,\gamma_J]=|\gamma_I|+|\gamma_J|.$ Indeed the function $\sigma(I,J)$ satisfies
$$
\sigma(I,J)=|I||J|-|I\cap J|-\sigma(I,J),
$$
and
\begin{eqnarray}
f(J,I) &=& (-1)^{\sigma(J,I)}\big(1-(-1)^{|I\cap J|}\big)\nonumber\\
    &=& (-1)^{|I||J|-|I\cap J|-\sigma(I,J)}\big(1-(-1)^{|I\cap J|}\big)\nonumber\\
    &=& (-1)^{|I||J|}(-1)^{\sigma(I,J)}\big((-1)^{|I\cap J|}-1\big)
    =-(-1)^{|I||J|}f(I,J). \nonumber
\end{eqnarray}
Hence $[\gamma_I,\gamma_J]=-(-1)^{|I||J|}[\gamma_J,\gamma_I]$ which shows that the relation (\ref{binary commutators}) is consistent with the symmetries of graded commutator. It is obvious that if the intersection of subsets $I,J$ contains an even number of elements then $f(I,J)=0$, and the graded commutator of $\gamma_I,\gamma_J$ is trivial. Particularly if at least one of two subsets $I,J$ is the empty set then $f(I,J)=0$. Thus any graded commutator (\ref{binary commutators}) containing $e$ is trivial.

As an example, consider the super Lie algebra $\frak C_2$. Its underlying vector space is 4-dimensional and $\frak C_2$ is generated by two even degree generators $e,\gamma_{12}$ and two odd degree generators $\gamma_1,\gamma_2.$ The non-trivial relations of this Lie superalgebra are given by
\begin{eqnarray}
[\gamma_1,\gamma_1]=[\gamma_2,\gamma_2]=2\,e,\;
[\gamma_1,\gamma_{12}]=2\,\gamma_2,\; [\gamma_2,\gamma_{12}]=-2\,\gamma_1.
\end{eqnarray}

Now we assume that $n=2m, m\geq 1$ is an even integer. The Lie superalgebra ${\frak C}_n$ has a matrix representation which can be described as follows. Fix $n=2$ and identify the generators $\gamma_1,\gamma_2$ with the Pauli matrices $\sigma_1,\sigma_2$, i.e.
\begin{equation}
\gamma_1=\left(
           \begin{array}{cc}
             0 & 1 \\
             1 & 0 \\
           \end{array}
         \right),\quad
\gamma_2=\left(
           \begin{array}{cc}
             0 & -i \\
             i & 0 \\
           \end{array}
         \right).
\label{Pauli matrices 1,2}
\end{equation}
Then $\gamma_{12}=\gamma_1\gamma_2=i\,\sigma_3$ where
$$
\sigma_3=\left(
           \begin{array}{cc}
             1 & 0 \\
             0 & -1 \\
           \end{array}
         \right).
$$
Let $S^2$ be the 2-dimensional complex super vector space $\mathbb C^2$ with the odd degree operators (\ref{Pauli matrices 1,2}), where the $\mathbb Z_2$-graded structure of $S^2$ is determined by $\sigma_3=i^{-1}\gamma_{12}$. Then $C_2\simeq \mbox{End}\,(S^2)$, and $S^2$ can be considered as a supermodule over the superalgebra $C_2$. Let $S^n=S^2\otimes S^2\otimes\ldots \otimes S^2 (m-\,\mbox{times})$. Then $S^n$ can be viewed as a supermodule over the $m$-fold tensor product of $C_2$, which can be identified with $C_n$ by identifying $\gamma_1,\gamma_2$ in the $j$th factor with $\gamma_{2j-1},\gamma_{2j}$ in $C_n$. This $C_n$-supermodule $S^n$ is called the supermodule of spinors \cite{Quillen-Mathai}. Hence we have the matrix representation for the Clifford algebra $C_n$, and this matrix representation or supermodule of spinors allows one to consider the supertrace, and it can be proved \cite{Quillen-Mathai} that
\begin{equation}
\mbox{Str}(\gamma_I)=\left\{
\begin{array}{ll}
0 & \text{if } I<{\mathcal N},\\
(2i)^m & \text{if } I={\mathcal N}.
\end{array} \right.
\label{super trace}
\end{equation}
Now we have the Lie superalgebra $\frak C_n$ with the graded commutator defined in (\ref{binary commutators}) and its matrix representation based on the supermodule of spinors. Hence we can construct a 3-Lie superalgebra by making use of graded ternary commutator (\ref{n super-to-n+1 super}). Applying the formula (\ref{n super-to-n+1 super}) we define the graded ternary commutator for any triple $\gamma_I,\gamma_J,\gamma_K$ of elements of basis for $\frak C_n$ by
\begin{eqnarray}
&& [\gamma_I,\gamma_J,\gamma_K] = \mbox{Str}(\gamma_I)\,[\gamma_J,\gamma_K]-(-1)^{|I||J|}\mbox{Str}(\gamma_J)\, [\gamma_I,\gamma_K]\nonumber\\
&&\qquad\qquad\qquad\qquad\qquad\qquad\quad + (-1)^{|K|(|I|+|J|)}\mbox{Str}(\gamma_K)\,[\gamma_I,\gamma_J],
\label{ternary graded commutator with super trace}
\end{eqnarray}
where the binary graded commutator at the right-hand side of this formula is defined by (\ref{binary commutators}). According to Theorem 2 the vector space spanned by $\gamma_I, I\subset {\mathcal N}$ and equipped with the ternary graded commutator (\ref{ternary graded commutator with super trace}) is the super 3-Lie algebra which will be denoted by ${\frak C}^{(3)}_n$. Making use of (\ref{binary commutators}) we can write the expression at the right-hand side of the above formula in the form
\begin{eqnarray}
&& [\gamma_I,\gamma_J,\gamma_K] =f(J,K) \mbox{Str}(\gamma_I)\,\gamma_{J\Delta K}-(-1)^{|I||J|}f(I,K)\mbox{Str}(\gamma_J)\, \gamma_{I\Delta K}\nonumber\\
&&\quad\qquad\qquad\qquad\qquad\qquad\qquad\quad + (-1)^{|K|(|I|+|J|)}f(I,J)\mbox{Str}(\gamma_K)\,\gamma_{I\Delta J}.\nonumber
\end{eqnarray}
From the formula for supertrace (\ref{super trace}) it follows immediately that the above graded ternary commutator is trivial if none of subsets $\gamma_i,\gamma_J,\gamma_K$ is equal to ${\mathcal N}$. Similarly this graded ternary commutator is also trivial if all three subsets $I,J,K$ are equal to $\mathcal N$, i.e. $I=J=K={\mathcal N}$, or two of them are equal to
 ${\mathcal N}$.
\begin{proposition}
The graded ternary commutators of the generators $\gamma_I,I\subseteq {\mathcal N}$ of the super 3-Lie algebra ${\frak C}^{(3)}_n$ are given by
\begin{equation}
[\gamma_I,\gamma_J,\gamma_K]=\left\{
\begin{array}{ll}
(2i)^mf(I,J)\gamma_{I\Delta J} & \text{if } I\neq {\mathcal N}, J\neq {\mathcal N}, K={\mathcal N},\\
0& \text{in all other cases }.
\end{array} \right.
\label{proposition}
\end{equation}
\end{proposition}
\begin{acknowledgement}
The authors is gratefully acknowledge the Estonian Science Foundation for financial support of this work under the Research Grant No. ETF9328. This research was also supported by institutional research funding IUT20-57
of the Estonian Ministry of Education and Research.
\end{acknowledgement}

\end{document}